\documentclass[number,citesort,seceqn,dvips]{arxbj}

\aid{0}
\volume{18}
\issue{1}
\pubyear{2012}
\firstpage{229}
\lastpage{251}
\doi{10.3150/10-BEJ326}

\makeatletter
\newcommand{\Ny}{N\rightarrow\infty}
\newcommand{\boldgreek}[1]{\bolds#1}
\newcommand{\Q}{\mathbb{P}}
\newcommand{\E}{\mathbb{E}}
\newcommand{\diag}{\operatorname{Diag}}
\newcommand{\matice}{{\boldgreek\Sigma}_2}
\newcommand{\X}{\mathbf{X}}
\newcommand{\Y}{\mathbf{Y}}
\newcommand{\cal}{\mathcal}
\newcommand{\mm}{\bolds\mu}
\newtheorem{LEM}{Lemma}[section]
\newremark{REM}{Remark}[section]
\newcommand{\eqref}[1]{(\ref{#1})}
\makeatother

\begin{document}
\begin{frontmatter}

\title{Nonparametric multivariate rank tests and their unbiasedness}
\runtitle{Nonparametric multivariate rank tests}

\begin{aug}
\author[a]{\fnms{Jana} \snm{Jure\v{c}kov\'{a}}\corref{}\thanksref{a}\ead[label=e1]{jurecko@karlin.mff.cuni.cz}} \and
\author[b]{\fnms{Jan} \snm{Kalina}\thanksref{b}\ead[label=e2]{kalina@euromise.cz}}
\runauthor{J. Jure\v{c}kov\'{a} and J. Kalina}
\pdfauthor{Jana Jureckova, Jan Kalina}
\address[a]{Department of Probability and Statistics, Charles
University in Prague,
Sokolovsk\'{a} 83, CZ-186 75 Prague~8, Czech Republic. \printead{e1}}
\address[b]{EUROMISE Center, Department of Medical Informatics,
Institute of Computer Science of the Academy of Sciences
of CR, v.v.i.,
Pod Vod\'{a}renskou v\v{e}\v{z}\'{i} 2, CZ-182 07 Prague 8,
Czech Republic.\\ \printead{e2}}
\end{aug}

\received{\smonth{11} \syear{2008}}
\revised{\smonth{9} \syear{2010}}

%
\begin{abstract}
Although unbiasedness is a basic property of a good test, many tests on
vector parameters or scalar parameters
against two-sided alternatives are not finite-sample unbiased. This was
already noticed by Sugiura
[\textit{Ann. Inst. Statist. Math.} \textbf{17} (1965) 261--263];
he found an alternative against which the Wilcoxon test is not
unbiased. The problem is even more serious in multivariate models. When
testing the hypothesis against an alternative which fits well with the
experiment, it should be verified whether the power of the test under
this alternative cannot be smaller than the significance level.
Surprisingly, this serious problem is not frequently considered in the
literature.

The present paper considers the two-sample multivariate testing
problem. We construct
several rank tests which are finite-sample unbiased against a broad
class of location/scale alternatives and are finite-sample
distribution-free under the hypothesis and alternatives. Each of them
is locally most powerful against a specific alternative of the Lehmann
type. Their powers against some alternatives are numerically compared
with each other and with other rank and classical tests. The question
of affine invariance of two-sample multivariate tests is also discussed.
\end{abstract}

%
\begin{keyword}
\kwd{affine invariance}
\kwd{contiguity}
\kwd{Kolmogorov--Smirnov test}
\kwd{Lehmann alternatives}
\kwd{Liu--Singh test}
\kwd{Psi test}
\kwd{Savage test}
\kwd{two-sample multivariate model}
\kwd{unbiasedness}
\kwd{Wilcoxon test}
\end{keyword}

\end{frontmatter}

\section{Introduction}\label{sec1}

\subsection{Two-sample multivariate tests}\label{sec1.1}
A frequent practical problem is that we have two data clouds of
$p$-dimensional observations with generally unknown distributions $F$
and $G,$ and
we wish to test the hypothesis that they both come from the same
distribution $F\equiv G,$ continuous but unknown.
Desirable properties of a test of such a hypothesis $\mathbf H$ are:
(i) being distribution-free under~$\mathbf H$; (ii) being affine
invariant with respect to changes of coordinate system; (iii) being
consistent against any fixed alternative; and (iv) being finite-sample
unbiased against a~broad class of alternatives of interest.
Unfortunately, a test satisfying all these conditions does not exist in
the multivariate setup.

Many authors have tried to attack this problem, emphasizing some of the
above properties. Their ideas were often concentrated either on some
geometric entities of the data clouds or on the affine invariance of
the testing problem.
Naturally, the ranks or the signed ranks of geometric entities of data
are invariant under many transformations
and provide a useful and simple tool for testing.
The papers extending methods based on ranks or other nonparametric
methods to the multivariate setup
use data depths,
Oja medians, multivariate sign functions and other tools. In this
context, we should mention the papers by Chaudhuri and Sengupta \cite
{ChaudhuriSengupta93},
Choi and Marden \cite{ChoiMarden97,ChoiMarden05}, Hallin and
Pandaveine \cite{HallinPandaveine02},
Hetmansperger \textit{et al.} \cite{HetOja98},
Liu \cite{Liu88,Liu90}, Liu and Singh \cite{LiuSingh93}, Oja
and Randles~\cite{OjaRandles04}, Oja \textit{et al.} \cite{OjaTienari97}, Puri
and Sen \cite{PuriSen71}, Randles and Peters \cite{RandlesPeters90},
Topchii \textit{et al.} \cite{Topchii03}, Tukey~\cite{Tukey75}, Zuo and He
\cite{ZuoHe06} and a recent excellent review by Oja \cite{Oja2010}.

Other authors have constructed various permutation tests: Bickel \cite
{Bickel69}, Brown \cite{Brown82}, Hall and Tajvidi \cite
{HallTajvidi02}, Neuhaus and Zhu \cite{NeuhausZhu99}, Oja \cite
{Oja97}, Wellner \cite{Wellner79}
and others. Tests based on distances between observations were
considered by Baringhaus and Franz \cite{BaringhausFranz04}, Friedman
and Rafsky \cite{FriedmanRafsky79}, Henze \cite{Henze88}, Maa \textit{et
al.} \cite{Maa96}, Rosenbaum \cite{Rosenbaum05} and Schilling \cite
{Schilling86}; the latter also
compared the simulated powers of his test with that of the
Kolmogorov--Smirnov two-sample test.

The proposed tests were typically consistent against distant
alternatives and some of them were affine invariant. The authors often
derived the asymptotic null distributions of the test criteria and
some derived the asymptotic powers under contiguous alternatives.
Many authors illustrated the powers on the simulated data, often
normal, and compared them with the power of the Hotelling $T^2$ test.
However, only in exceptional cases did they check whether the test was
unbiased against alternatives of interest.

\subsection{Unbiased tests}\label{sec1.2}
Let $\Phi$ be a test of hypothesis $\{\mathbf H\dvtx $ distribution $F$ of
random vector $\mathbf X$ belongs to the set $\mathcal H\}$
against the alternative $\{\mathbf K\dvtx $ distribution of $\mathbf X$
belongs to the set $\mathcal K\}.$ Consider the tests of size $\alpha
,  0<\alpha<1,$ where $\alpha$ is the chosen significance level,
that is, the tests satisfying $\sup_{F\in\mathcal H}\E_F[\Phi
(\mathbf X)]\leq\alpha.$ The test $\Phi$ is unbiased if it satisfies
\[
\sup_{F\in\mathcal H}\E_F[\Phi(\mathbf X)]\leq\alpha  \quad \mbox{and} \quad
\inf_{F\in\mathcal K}\E_F[\Phi(\mathbf X)]\geq\alpha.
\]
This is a natural property of a test; it means that the power of a test
should not be smaller than the permitted error of the first kind. If
the test rejects the hypothesis with a probability less than
$\alpha$ under the alternative of interest, then we can hardly
recommend the test to the experimenter. Note that if there exists a
uniformly most powerful test, then it is always unbiased. If the
optimal test of size $\alpha$ does not exist because the family of
$\alpha$-tests is too broad, then we should restrict ourselves to a
pertinent subfamily of tests, and the family of unbiased tests of size
$\alpha$ is the most natural subfamily. We refer to Lehmann's
monograph \cite{Lehmann97} for an excellent account of unbiased tests.

Many tests criteria have asymptotic normal distributions under the
hypothesis as well as under the local alternatives -- these are \textit
{asymptotically} locally unbiased.
However, the practice always works with a finite number of
observations. The asymptotic distribution only approximates well the
central part of the finite-sample distribution; elsewhere, it can
stretch the truth and sometimes is only valid for a huge number of observations.
To calculate the finite-sample power of a test is sometimes difficult;
in any case, as a first step, we should be sure that the test is
unbiased against the alternatives under consideration, at least locally
in a neighborhood of the hypothesis. Unfortunately, many authors have
not specified the alternatives against which their tests are (locally)
unbiased. The alternative is typically more important for an
experimenter than the hypothesis because it describes his/her
scientific conjecture.
Some papers, for example, \cite{Jur2002,Xavier03,Sugiura65,Sugiura06}, show that the tests are not automatically
finite-sample unbiased. While the univariate two-sample Wilcoxon test,
for example, is always unbiased against one-sided alternatives, it is
generally not unbiased against two-sided alternatives, even not with
equal sample sizes (see \cite{Sugiura65,Sugiura06}). The test is locally unbiased against two-sample
alternatives only under some conditions on the hypothetical
distribution of observations (e.g., when it is symmetric). Amrhein
\cite{Amrhein95} demonstrated the same phenomenon for the one-sample
Wilcoxon test. Hence, the finite-sample unbiasedness of some tests
cited above, and others described in the literature, is still an open question.

To illustrate this problem more precisely, consider a random vector
${\mathbf X}=(X_1,\ldots ,X_n)$ with
distribution function $F(\mathbf x,{\boldgreek\theta}),  {\boldgreek
\theta}\in{\boldgreek\Theta}\subset\mathbb{R}^p$, and density
$f(\mathbf x,{\boldgreek\theta})$ (not necessarily Lebesgue), which
has bounded third derivatives in components of $\boldgreek\theta$ in
a neighborhood of~${\boldgreek\theta}_0$ and a positive definite
Fisher information matrix. We wish to test ${\mathbf H}_0\dvtx   {\boldgreek\theta}={\boldgreek\theta}_0$ against the alternative
${\mathbf K}\dvtx   {\boldgreek\theta}\neq{\boldgreek\theta}_0$ using
the test $\Phi$ of size $\alpha,$ that is, ${\E}_{\theta_0}[{\Phi
}({\mathbf X})]=\alpha.$
We then have the following expansion of the power function of $\Phi$
around $\boldgreek\theta_0$ (see \cite{Xavier03}):
%
\begin{eqnarray} \label{expansion}
 {\E}_{\theta_0}\Phi({\mathbf X})&=&\alpha+({\boldgreek\theta
}-{\boldgreek\theta}_0)^{\top}{\E}_{\theta_0} \biggl\{\Phi
({\mathbf X})
\frac{ (\dot{\mathbf f}({\mathbf X},{\boldgreek\theta
}_0) )}{f({\mathbf X},{\boldgreek\theta}_0)} \biggr\}\nonumber
\\[-8pt]
\\[-8pt]
&&{}+\frac12({\boldgreek\theta}-{\boldgreek\theta}_0)^{\top}{\E
}_{\theta_0} \biggl\{\Phi({\mathbf X})\frac{ [\ddot{\mathbf
f}({\mathbf X},{\boldgreek\theta}_0) ]}{f({\mathbf
X},{\boldgreek\theta}_0)} \biggr\}({\boldgreek\theta}-{\boldgreek
\theta}_0)+
{\mathcal O} (\|{\boldgreek\theta}-{\boldgreek\theta}_0\|
^3 ),
\nonumber
\end{eqnarray}
where
\[
 (\dot{\mathbf f}({\mathbf x},{\boldgreek\theta}) )=
\biggl(\frac{\partial f({\mathbf x},{\boldgreek\theta})}{\partial\theta
_1},\ldots ,\frac{\partial f({\mathbf x},{\boldgreek\theta
})}{\partial\theta_p} \biggr)^{\top}, \qquad
 [\ddot{\mathbf f}({\mathbf x},{\boldgreek\theta})
]= \biggl[\frac{\partial^2f({\mathbf x},{\boldgreek\theta
})}{\partial\theta_j\,\partial\theta_k} \biggr]_{j,k=1}^p.
\]
The test $\Phi$ is locally unbiased if the second term on the
right-hand side of (\ref{expansion}) is nonnegative.
If $\theta$ is a scalar parameter and we consider the one-sided
alternative \mbox{$\mathbf K\dvtx  \theta>\theta_0,$} then there always exists an
unbiased test. However, the alternative for a vector $\boldgreek\theta
$ is only two-sided and the local unbiasedness of $\Phi$
is guaranteed only when
%
\begin{equation}\label{condition}
{\E}_{\theta_0} \biggl\{\Phi({\mathbf X})
\frac{ (\dot{\mathbf f}({\mathbf X},{\boldgreek\theta
}_0) )}{f({\mathbf X},{\boldgreek\theta}_0)} \biggr\}={\mathbf0}.
\end{equation}
However, (\ref{condition}) is generally true only for $f$ satisfying
special conditions, which cannot easily be verified for unknown $f.$
If the test $\Phi$ does not satisfy (\ref{condition}), then the
second term in
(\ref{expansion}) can be negative for some ${\boldgreek\theta}$ and
hence the power of $\Phi$ can be less than~$\alpha.$
We refer to Grose and King \cite{GroseKing91}, who imposed condition
(\ref{condition}) when constructing a~locally unbiased two-sided
version of the Durbin--Watson test.

\subsection{Outline of the paper}\label{sec1.3}

We shall propose three classes of multivariate two-sample tests, based
on the ranks of suitable distances of multivariate observations. One
test is based on the ranks of distances of observations from the
origin, while the others are based on the ranks of their interpoint
distances. The natural alternatives state that either the distances of
the second sample from the origin are stochastically larger than those
of the first sample, or that the distances of the $\mathbf Y$'s from
the $\mathbf X$'s are stochastically larger than those of the~$\mathbf
X$'s from each other. The proposed tests are unbiased because our
natural alternatives are one-sided (in the distances).
Moreover, the proposed rank tests are 
distribution-free under the hypothesis as well as under
alternatives of the Lehmann type, and they are consistent against
general alternatives (properties (i), (iii) and (iv)). The
distribution-free property is important because we need not determine the
distribution of distances when performing the test.

The tests are described in Section \ref{sec2}, which starts with
some invariance considerations (cf. desired property (ii) of the test).
It is shown that the proposed tests based on the ranks of distances, as
well as the Liu--Singh \cite{LiuSingh93} tests based on the ranks of
depths, are distribution-free
against some monotone alternatives of the Lehmann type with respect to
which they are finite-sample unbiased. 
Section \ref{sec3} describes the contiguity of these alternatives
with respect to the hypothesis, which enables us to derive the local
asymptotic powers of the tests. The powers of tests are compared
numerically under finite $N,$ as well as asymptotically under $\Ny.$
The proposed tests are also compared with the tests of Liu and Singh,
and of
 \cite{HetOja98}, using a
reference to numerical results of   \cite{ZuoHe06}. In
Section~\ref{sec4}
we compare the empirical powers of the proposed tests with those of the
Hotelling test under the bivariate normal and bivariate Cauchy
distributions. The contiguity of the Lehmann-type alternatives is
proved in the \hyperref[appm]{Appendix}.

\section{Multivariate two-sample rank tests}\label{section2}\label{sec2}
\subsection{Remarks to affine invariance}\label{sec2.1}

Consider two independent samples ${\mathcal X}=({\mathbf X}_1,\ldots
,{\mathbf X}_m)$ and ${\mathcal Y}=({\mathbf Y}_1,\ldots ,{\mathbf
Y}_n)$ from two $p$-variate populations with
continuous distribution functions $F^{(p)}$ and $G^{(p)},$
respectively, with respective means and dispersion matrices $\boldgreek
{\mu}_1,\boldgreek{\mu}_2$ and $\boldgreek{\Sigma}_1,\boldgreek
{\Sigma}_2.$ The problem is to test the hypothesis ${\mathbf H}_0\dvtx
F^{(p)}\equiv G^{(p)} \mbox{ (along with }  \boldgreek{\mu
}_1=\boldgreek{\mu}_2,  \boldgreek{\Sigma}_1=\boldgreek{\Sigma
}_2)$ against an alternative $\mathbf{H}_1,$ where either $(\boldgreek
{\mu}_1,\boldgreek{\Sigma}_1)\neq(\boldgreek{\mu}_2,\boldgreek
{\Sigma}_2)$ or where $F^{(p)}$ and $G^{(p)}$ are not of the same
functional form.\vadjust{\goodbreak}
We denote by $(\mathbf{Z}_1,\ldots ,\mathbf{Z}_{N})$ the pooled
sample with $\mathbf{Z}_i=\mathbf{X}_i,  i=1,\ldots ,m,$ and
$\mathbf{Z}_{m+j}=\mathbf{Y}_j,  j=1,\ldots ,n,  m+n=N.$
The hypothesis and the alternative are invariant under affine transformations:
%
\begin{equation}\label{affine}
\mathcal{G}\dvtx  \{\mathbf{Z}\rightarrow\mathbf{a} +\mathbf{BZ}\}  \mbox{ with }\mathbf a\in\mathbb{R}^p  \mbox{and } \mathbf{B}
\mbox{ a nonsingular } p\times p  \mbox{ matrix}.
\end{equation}
More precisely, the hypothesis and alternative remain true even after
the transformation $g\in\mathcal G$ of the data, and we are looking
for invariant tests whose criteria are invariant with respect to $g\in
\mathcal G.$ The invariant tests
depend on the data only by means of a \textit{maximal invariant} of
$\mathcal{G}$  \cite{Lehmann97}. Obenchain \cite
{Obenchain71} showed that
the maximal invariant with respect to~$\mathcal{G}$~is\looseness=-1
%
\begin{eqnarray} \label{21b}
&&\mathbf{T}(\mathbf{Z}_1,\ldots ,\mathbf{Z}_N)= [(\mathbf
{Z}_i-\bar{\mathbf Z}_N)^{\top}\mathbf{V}_N^{-1}(\mathbf{Z}_j-\bar
{\mathbf Z}_N) ]_{i,j=1}^N,\nonumber
\\[-8pt]
\\[-8pt]
&& \quad \mbox{where }
\bar{\mathbf Z}_N= \frac1N\sum_{i=1}^N \mathbf{Z}_i,  \mathbf
{V}_N=\sum_{i=1}^N (\mathbf{Z}_i-\bar{\mathbf Z}_N)(\mathbf
{Z}_i-\bar{\mathbf Z}_N)^{\top}.
\nonumber
\end{eqnarray}\looseness=0
Then, $\mathbf{T}(\mathbf{Z}_1,\ldots ,\mathbf{Z}_N)$ is the
projection matrix associated with the space spanned by the columns of
the matrix $ [\mathbf{Z}_1-\bar{\mathbf Z}_N,\ldots ,\mathbf
{Z}_N-\bar{\mathbf Z}_N ].$ In particular,
under $\mathbf a\equiv\mathbf0,$ the maximal invariant of the group
%
\begin{eqnarray}\label{21}
&&\mathcal{G}_0\dvtx  \{\mathbf{Z}\rightarrow\mathbf{BZ}\}  \quad
\mbox{is equal to} \quad \mathbf{T}_0(\mathbf{Z}_1,\ldots ,\mathbf
{Z}_N)= [\mathbf{Z}_i^{\top}(\mathbf{V}_N^{0})^{-1}
\mathbf{Z}_j ]_{i,j=1}^N,\nonumber
\\[-8pt]
\\[-8pt] && \quad     \mbox{where }  \mathbf{V}_N^0=\sum
_{i=1}^N \mathbf{Z}_i\mathbf{Z}_i^{\top}.
\nonumber
\end{eqnarray}
Moreover, one of the maximal invariants with respect to the \textit
{group of shifts in location}
%
\begin{equation}\label{G*}
\mathcal{G}_{1}\dvtx   \mathbf Z\longrightarrow\mathbf Z+\mathbf a,
\qquad
 \mathbf a\in\mathbb{R}^p,
\end{equation}
is $\mathbf T_1(\mathbf{Z}_1,\ldots ,\mathbf{Z}_N)=(\mathbf
Z_2-\mathbf{Z}_1,\ldots ,\mathbf{Z}_N-\mathbf{Z}_1).$

The well-known two-sample Hotelling $T^2$ test is based on the criterion
%
\begin{equation}\label{Hotelling}
\mathcal T_{mn}^2=(\bar{\mathbf X}_m-\bar{\mathbf Y}_n)^{\top
}\mathbf V_N^{-1}(\bar{\mathbf X}_m-\bar{\mathbf Y}_n).
\end{equation}
The test is invariant with respect to $\mathcal G$ and is optimal
unbiased against two-sample normal alternatives with
$\boldgreek{\mu}_1\neq\boldgreek{\mu}_2$ and $\boldgreek\Sigma
_1=\boldgreek\Sigma_2.$
Its asymptotic null distribution, when both sample sizes $m,n$ tend to
infinity, does not depend on the normality. If $m,n\rightarrow\infty$
and $\frac{m}{n}\rightarrow1,$ then the asymptotic distribution of
$\mathcal T_{mn}^2$ does not change even when $\boldgreek\Sigma_1\neq
\boldgreek\Sigma_2,$ but only in this case (see \cite{Lehmann97}).
Its finite sample unbiasedness is not guaranteed under a~nonnormal
underlying distribution.

If we wish to construct a nonparametric two-sample test which is
distribution-free and affine invariant with respect to $\mathcal G,$ we
expect it to be based on the ranks of some components of $\mathbf T$ in
(\ref{21b}) or on the relevant empirical Mahalanobis distances of
points~$\mathbf Z_i, \mathbf Z_j,  1\leq i,j\leq N.$ The ranks of
distances are invariant with respect to continuous increasing functions
of the distances; however, in our case, the data themselves are
transformed,
rather than their distances. The proper form of the rank test criterion
based on the Mahalanobis distances\vadjust{\goodbreak} and its unbiasedness against
alternatives of interest is the subject of a forthcoming study. The
rank tests considered in the present paper are easier but invariant
only with respect to $\mathcal G_1,$ not to the change of the origin.
On the other hand, the proposed tests enjoy the properties (i), (iii)
and (iv) mentioned above.

\subsection{Liu and Singh rank sum test}\label{sec2.2}

An interesting test of Wilcoxon type, based on the ranks of depths of
the data, was proposed by Liu and Singh \cite{LiuSingh93}. Being of
Wilcoxon type, this test
is locally most powerful against some alternatives of the Lehmann type.
Its asymptotic distributions under the hypothesis and under general
alternative distributions $F,G$ of depths was derived by Zuo and He
\cite{ZuoHe06}.

Let $D(\mathbf y;H)$ denote a depth function of a distribution $H$
evaluated at the point $\mathbf y\in\mathbb{R}^p.$ Liu and Singh
\cite{LiuSingh93} considered a parameter, called a \textit{quality
index}, defined as
\begin{eqnarray*}
  Q\bigl(F^{(p)},G^{(p)}\bigr)&=&\int R\bigl({\mathbf y};F^{(p)}\bigr)\,\mathrm{d}G^{(p)}({\mathbf y})\\
    &=&{\Q} \bigl\{ D\bigl({\mathbf X};F^{(p)}\bigr)\leq
D\bigl({\mathbf Y};F^{(p)}\bigr) |{\mathbf X}\sim F^{(p)}, {\mathbf Y}\sim
G^{(p)} \bigr\},
\end{eqnarray*}
where $R({\mathbf y};F^{(p)})={\Q}_F (D({\mathbf X};F^{(p)})\leq
D({\mathbf y};F^{(p)} ),  {\mathbf y}\in\mathbb{R}^p,$ and
showed that if $D({\mathbf X};F^{(p)})$ has a continuous distribution,
then $Q(F^{(p)}, F^{(p)})=\frac12.$ They then tested the hypothesis
$Q(F^{(p)},G^{(p)})=\frac12$ against the alternative
$Q(F^{(p)},G^{(p)})\neq\frac12$ using the Wilcoxon-type criterion
based on the empirical distribution functions $F_m$ and $G_n$ of
samples of sizes $m$ and $n$ respectively:
\[
Q(F_m,G_n)=\int R({\mathbf y};F_m)\,\mathrm{d}G_n({\mathbf y})=\frac1n\sum
_{j=1}^n R({\mathbf Y}_j;F_m).
\]
If the distribution of depths is symmetric under $F^{(p)}\equiv
G^{(p)},$ then the test rejecting provided $|Q(F_m,G_n)-\frac12|\geq
C_{\alpha/2}$ is locally unbiased against $Q(F^{(p)},G^{(p)})\neq
\frac12.$ Under a~general distribution of depths, only the one-sided
test with the critical region
\[
Q(F_m,G_n)-\tfrac12 >C_{\alpha}
\]
is unbiased against the one-sided alternative $Q(F^{(p)},G^{(p)})>\frac
12;$ however, this alternative, one-sided in depths, has a difficult
interpretation with respect to the distributions $F^{(p)}$ and
$G^{(p)}$ of original observations ${\mathbf X}$ and ${\mathbf Y}$,
respectively.
Generally, the test is not finite-sample unbiased against $F\neq G,$
not even locally. The unbiasedness can be guaranteed only in some
cases, for instance, if the hypothetical distribution of depths is symmetric.

\subsection{Rank tests based on distances of observations}\label{sec2.3}

We shall test the hypothesis of equality of distributions of two
samples against alternatives that some distances are greater than
others; because such alternatives are one-sided, they make the tests unbiased.

Choose a distance $L=L(\cdot,\cdot)$ in ${\mathbb R}^p$ taking
nonnegative real values.
Let
\[
{\mathcal Z}=({\mathbf Z}_1,\ldots ,{\mathbf Z}_N)=({\mathbf
X}_1,\ldots ,{\mathbf X}_m, {\mathbf Y}_1,\ldots ,{\mathbf Y}_n)
\]
denote the pooled sample, where $N=m+n,$ and consider the matrix of
distances $\mathbb L_N=[\ell_{ik}]_{i,k=1}^N,$ where $\ell
_{ik}=L(\mathbf Z_i,\mathbf Z_k).$ We can construct simple rank tests
based on $\mathbb L_N$ in three different ways:
\begin{enumerate}[(iii)]
\item[(i)] \textit{Simple rank test, but not invariant with respect
to $\mathcal G$ or $\mathcal G_1.$}
Consider the vector
\[
(\tilde{\ell}_1,\ldots ,\tilde{\ell}_N), \qquad   \tilde{\ell
}_k=L(\mathbf0,\mathbf Z_k), \qquad   k=1,\ldots ,N,
\]
of distances of observations from the origin. The vector $(\tilde{\ell
}_1,\ldots ,\tilde{\ell}_m)$ is then a~random sample from
a population with a distribution function $F$ (say), while $(\tilde
{\ell}_{m+1},\ldots ,\tilde{\ell}_N)$ is a random sample from
a population with a distribution function~$G.$ Assume that the
distribution functions $F$ and $G$ are absolutely continuous. Under
hypothesis ${\mathbf H}_0\dvtx   F^{(p)}\equiv G^{(p)},$ the distribution
functions $F$ and $G$ coincide, that is, they satisfy the hypothesis
$\widetilde{\mathbf H}_0\dvtx   F\equiv G$ which states that $\{\tilde
{\ell}_{k},  k=1,\ldots ,N \}$ satisfy the hypothesis of randomness.
If $\widetilde{\mathbf H}_0$ is not true, then ${\mathbf H}_0$ is not
true either.

Denote by $(\tilde{R}_{1},\ldots ,\tilde{R}_N)$ the respective ranks
of $\{\tilde{\ell}_{k},  k=1,\ldots ,N\}.$
Under the hypothesis, the vector of ranks has the uniform distribution
on the set of permutations of the numbers $1,\ldots ,N.$
Because $\{\tilde{\ell}_{k},   k=1,\ldots ,m\}$ and $\{\tilde{\ell
}_{k},  k=m+1,\ldots ,N\}$ are random samples,
under the hypothesis, as well as under the alternatives, every
two-sample rank test will depend only on the ordered ranks
$R^{(m+1)}<\cdots <R^{(N)}$ of the second sample. However, although
invariant with respect to increasing continuous functions of $(\tilde
{\ell}_1,\ldots ,\tilde{\ell}_N),$ such a test would not be
invariant with respect to the groups of transformations (\ref{affine})
or (\ref{21}),
even if $\tilde{\ell}_k=\|\mathbf Z_k\|$ is the Euclidean distance.

The linear rank test is based on the linear rank statistic
%
\begin{equation}
\label{1}
S_{N}=N^{-1/2}\sum_{k=m+1}^N a_N(R_{ik})
\end{equation}
with the scores $a_N(1),\ldots ,a_N(N)$
generated by a nondecreasing score function
$\varphi: (0,1)\mapsto\mathbb R$ in either of the following two ways:
%
\begin{equation}\label{exact}
a_N(k)= {\E} (\varphi(U_{N:k}) ), \qquad   k=1,\ldots ,N,
\end{equation}
or
%
\begin{equation}\label{approx}
a_N(k)=\varphi \biggl(\frac{U_{N:k}}{N+1} \biggr), \qquad   k=1,\ldots ,N,
\end{equation}
with $U_{N:1}\leq\cdots\leq U_{N:N}$ being the order statistics of
the sample from the uniform ${\mathbf R}(0,1)$ distribution. The test
based on (\ref{1}) is distribution-free, that is, the null
distribution of $S_{N}$ does not depend on the unknown $F\equiv G$
under $\widetilde{\mathbf H}_0.$ Its asymptotic properties follow from
\cite{HajekSidak67} or \cite{HajekSidakSen99}.

\item[(ii)] \textit{Conditional rank test, invariant with respect to
$\mathcal G_1.$}
Assuming that $m>p,$ choose a suitable basis $(\mathbf X_{i_1},\ldots
,\mathbf X_{i_p})=\mathbb X_p$ of $\{\mathbf{X}_i,  1\leq i\leq m\};$
the choice of basis $\mathbb X_{p}$
can follow various aspects. Consider the set of $(m+n-p)\times p$ distances
\[
\{\ell_{i_j,k}^*=L({\mathbf X}_{i_j},{\mathbf Z}_k),   k=1,\ldots ,N,
 k\neq i_1,\ldots ,i_p\},  \qquad  j=1,\ldots ,p.
\]
Then, for a fixed $i_j,  1\leq j\leq p,$ and conditionally given
$\mathbb X_p,$ the vector $\{\ell_{i_j,k}^*,   k=1,\ldots ,m,  k\neq
i_1,\ldots ,i_p\}$ is a random sample from a population with a
distribution function $F(z|\mathbb{X}_p)=F$ (say), while $\{\ell
_{i_j,k}^*,  k=m+1,\ldots ,N\}$
is a random sample from a population with a distribution function
$G(z|\mathbb{X}_p)=G.$ Assume that the distribution functions $F$ and
$G$ are absolutely continuous.
Let
\[
{\mathcal R}_{i_j}=(R_{i_j,k},  k=1,\ldots ,N,  k\neq i_1,\ldots ,i_p)
\]
denote the ranks of $\ell_{i_jk}^*,  k=1,\ldots ,N,  k\neq i_j,
 \forall j=1,\ldots ,p.$ Every two-sample rank test will depend only on
the ordered ranks $R_{i_j}^{(m+1)}<\cdots  <R_{i_j}^{(N)}$ of the second
sample. In particular, if $L(\mathbf X_{i_j}, \mathbf Z_k)=\|\mathbf
X_{i_j}-\mathbf Z_k\|,  k=1,\ldots ,N,  k\neq i_1,\ldots ,i_p,$ where~$\|\cdot\|$ is the Euclidean distance, then the test based on their
ranks will be invariant with respect to $\mathcal G_1$ in (\ref{G*}),
but not with respect to $\mathcal G,  \mathcal G_0.$

Similarly as in {(i)}, the linear (conditional) rank test is
based on the linear rank statistic
\[
S^*_{i_j,N}=N^{-1/2}\sum_{k=m+1}^N a_N(R_{i_j,k})
\]
with the scores $a_N(1),\ldots ,a_N(N-p)$ generated by a nondecreasing
score function~$\varphi$ as in either (\ref{exact}) or (\ref
{approx}). The criteria $S^*_{i_j,N}$ are equally distributed for
$j=1,\ldots ,p$ under the hypothesis and under the alternatives, and
are conditionally independent given $\mathbb X_p.$ Using only a single
$S^*_{i_j,N}$ would be a loss of information, so we look for a
convenient combination of $S^*_{i_1,N},\ldots ,S^*_{i_p,N}.$
Every convenient homogeneous combination of $S^*_{i_1,N},\ldots
,S^*_{i_p,N}$ leads to a rank test,
conditional under given $\mathbb X_p,$ which is distribution-free under
the hypothesis. The problem would be to find its null distribution, and
thus the critical values, under a finite~$N.$
The test based on a single $S^*_{i_j,N}$ is a standard rank test, for
example, Wilcoxon, conditionally given $\mathbb X_p,$ and is thus easy
to perform. When we look for a similarly simple test based on a
combination of $S^*_{i_jN},  1\leq j\leq p,$ it seems that the simplest
possibility is a randomization of $S^*_{i_1,N},\ldots ,S^*_{i_p,N},$
leading to the following criterion~$\widetilde{S}^{(N)}$:
%
\begin{equation}
\label{3a}
{\Q}\bigl(\widetilde{S}^{(N)}=S^*_{i_j,N}\bigr)=\frac1p, \qquad    j=1,\ldots ,p,
\end{equation}
where the randomization in (\ref{3a}) is independent of the set of
observations ${\mathcal Z}.$
The following identity is true for any $C$:
\[
{\Q} \bigl(\widetilde{S}^{(N)}>C \bigr)=\frac1p\sum_{j=1}^p{\Q
} (S^*_{i_j,N}>C )
\]
and the test rejects $\widetilde{\mathbf H}_0$ for $\alpha\in(0,1)$
if $\widetilde{S}^{(N)}>C_{\alpha};$ eventually, it rejects with probability
$\gamma\in(0,1)$ if $S^{(N)}=C_{\alpha},$ where
\[
{\Q}_{\widetilde{\mathbf H}_0} \bigl(\widetilde{S}^{(N)}>C_{\alpha
} \bigr)+\gamma{\Q}_{\widetilde{\mathbf H}_0}\bigl (\widetilde
{S}^{(N)}=C_{\alpha} \bigr)=\alpha.
\]

\item[(iii)] \textit{Randomized rank test, invariant with respect to
$\mathcal G_1.$}
Similarly, for every fixed $i$ and under fixed $\mathbf{X}_i,  1\leq
i\leq m,$ we can consider the distances $\{\ell_{ik}^*=L({\mathbf
X}_i,{\mathbf Z}_k),   k=1,\ldots ,N,  k\neq i\}.$ Then, conditionally
given $\mathbf X_i,$ the vector $\{\ell_{ik}^*,   k=1,\ldots ,m,
 k\neq i\}$ is a random sample from a population with a distribution
function $F(z|\mathbf{X}_i)=F$ (say), while $\{\ell_{ik}^*,  k=m+1,\ldots ,N\}$
is a random sample from a population with a distribution function
$G(z|\mathbf{X}_i)=G.$
Assuming that the distribution functions $F$ and $G$ are absolutely
continuous, we work with
\[
{\mathcal R}_i=(R_{i1},\ldots ,R_{i,i-1},R_{i,i+1},\ldots ,R_{iN}),
\]
the ranks of $\ell_{ik}^*,  k=1,\ldots ,N,  k\neq i.$

The linear (conditional) rank test is based on the linear rank statistic
%
\begin{equation}
\label{10}
S_{iN}=N^{-1/2}\sum_{k=m+1}^N a_N(R_{ik})
\end{equation}
with the scores $a_N(1),\ldots ,a_N(N-1).$

The criteria $S_{iN}$ are equally distributed for $i=1,\ldots ,m$ under
the hypothesis and under the alternatives, although not independent.
We look for a convenient combination of $S_{1N},\ldots ,S_{mN}.$
Again, a randomization of $S_{1N},\ldots ,S_{mN}$ keeps the simple
structure of the test and is thus easy to perform. It leads to the
following criterion, $S^{(N)}$:
%
\begin{equation}
\label{3}
{\Q}\bigl(S^{(N)}=S_{iN}\bigr)=\frac1m, \qquad    i=1,\ldots ,m,
\end{equation}
where the randomization in (\ref{3}) is independent of the set of
observations ${\mathcal Z}.$
Again, for any $C,$
\[
{\Q}\bigl (S^{(N)}>C \bigr)=\frac1m\sum_{i=1}^m{\Q}
(S_{iN}>C ),
\]
and the test rejects $\widetilde{\mathbf H}_0$ for $\alpha\in(0,1)$
if $S^{(N)}>C_{\alpha};$ eventually, it rejects with probability
$\gamma\in(0,1)$ if $S^{(N)}=C_{\alpha}.$
Again, with the Euclidean distance, the test will be invariant with
respect to $\mathcal G_1$ in (\ref{G*}), but not with respect to
$\mathcal G,  \mathcal G_0.$
\end{enumerate}

\begin{REM}\label{Mahal} The Mahalanobis distances
%
\begin{eqnarray}\label{Mahalanobis}
&\displaystyle\mathbf Z_k^{\top}(\mathbf{V}_N^0)^{-1}\mathbf Z_k,  \qquad  k=1,\ldots
,N,&\\[2pt]
\label{Mahala}
&\displaystyle(\mathbf X_i-\mathbf Z_k)^{\top}\mathbf{V}_N^{-1}(\mathbf
X_i-\mathbf Z_k)   \quad \mbox{or} \quad   (\mathbf X_i-\mathbf Z_k)^{\top}(\mathbf{V}_N^0)^{-1}(\mathbf
X_i-\mathbf Z_k),  \qquad   k\neq i,&
\end{eqnarray}
are not independent, but under $\mathbf H_0$, they have exchangeable
distributions; hence, under~$\mathbf H_0$, the distribution of
their
ranks is independent of the distribution of observations (is distribution-free).
Moreover, (\ref{Mahala})
are invariant with respect to $\mathcal G$ and $\mathcal G_0,$ while
(\ref{Mahalanobis}) are invariant only with respect to $\mathcal G_0.$
The invariant tests based on the ranks of~(\ref{Mahalanobis}) or~(\ref
{Mahala})
will be the subject of a further study. Their structure is more complex
than that of tests based on simple distances.\vspace*{-2pt}
\end{REM}
%
\section{Structure of the rank tests}\vspace*{-2pt}\label{section2a}\label{sec3}
Let $\mathbf X=(X_1,\ldots ,X_m)$ and $\mathbf Y=(Y_1,\ldots ,Y_n)$ be
two independent samples from distributions $F$ and $G$, respectively.
Consider the rank test with the criterion $S_N=N^{-1/2}\sum
_{k=m+1}^Na_N(R_i),$ where $R_1,\ldots ,R_N$ are the ranks of the
pooled sample $\mathbf Z=(X_1,\ldots ,X_m,Y_1,\ldots ,Y_n).$ The values
$X_i,  Y_j$ are, for example, the distances of multivariate
observations, either from a fixed point or the interpoint distances
considered conditionally given the original component. 
We want to test the hypothesis $\mathbf H_0\dvtx   F\equiv G$ against
a~general alternative with the $(m+n)$-dimensional distribution function
of the form
%
\begin{equation}\label{alternatives}
{\mathbf K}\dvtx   \prod_{k=1}^m G_{\Delta}^{(1)}(z_k)\prod
_{k=m+1}^NG_{\Delta}^{(2)}(z_k).
\end{equation}
Lehmann \cite{Lehmann53} showed that the Wilcoxon two-sample test with
the score-generating function $\varphi(u)=2u-1,  0\leq u\leq1,$ is
the locally most powerful rank test of ${\mathbf H}_0$ against the
class of alternatives
with
%
\begin{eqnarray}\label{51}
 G_{\Delta}^{(1)}(z)&=&F(z)   \quad \mbox{and} \quad   G_{\Delta
}^{(2)}(z)=G_{\Delta}(z),\\
  \label{Lehmann1} 
 G_{\Delta}(z)&=&
\cases{\displaystyle  (1-\Delta) F(z)+\Delta F^2(z), & \quad  $z\geq0$,\cr\displaystyle
0, & \quad  $z<0$
}
\end{eqnarray}
with $0< \Delta<1.$ The $Y$'s are then stochastically larger than the
$X$'s and
$F(z)-G_{\Delta}(z) \equiv\Delta\cdot F(z)(1-F(z));$ hence, the
Kolmogorov distance of $F$ and $G_{\Delta}$ is
\[
d_K(F,G_{\Delta})=\Delta\cdot\sup_{z\geq0}
\bigl[F(z)\bigl(1-F(z)\bigr) \bigr]= \frac{\Delta}{4}
\]
and the point of maximum is $z=F^{-1}(\frac12).$

Gibbons \cite{Gibbons64}
proved that the \textit{Psi test} with the scores
\[
a_N(i)=\sum_{j=0}^{i-1}\frac{1}{N-j}-\sum_{j=0}^{N-i}\frac
{1}{N-j}, \qquad   i=1,\ldots ,N,
\]
is the locally most powerful rank test of ${\mathbf H}_0$
against the alternative (\ref{alternatives}) with
%
\begin{eqnarray}\label{G}
 G_{\Delta}^{(1)}(z)&=&1 - \bigl(1 - F(z)\bigr)^{1+\Delta}, \nonumber
 \\[-8pt]
 \\[-8pt]
 G_{\Delta}^{(2)}(z)&=& (F(z))^{1+\Delta},  \qquad  \Delta>0,\  z\geq
0.
\nonumber\vadjust{\goodbreak}
\end{eqnarray}
Obviously,
\begin{eqnarray*}
G_{\Delta}^{(1)}(z)&\geq& F(z)\geq G_{\Delta}^{(2)}(z)\\[-1pt]
&=&
(F(z))^{1+\Delta}
 \qquad  \forall z\geq0  \mbox{ and } \Delta\geq0,\vspace*{-1pt}
\end{eqnarray*}
hence $G_{\Delta}^{(1)}$ is stochastically smaller than $G_{\Delta
}^{(2)}$ for $\Delta>0.$
The Kolmogorov distance of~$G_{\Delta}^{(1)}$ and $G_{\Delta}^{(2)}$ is
\begin{eqnarray*}
d_K \bigl(G_{\Delta}^{(1)},G_{\Delta}^{(2)} \bigr)&=&\sup_{z\geq
0} \bigl[1-(F(z))^{1+\Delta}-
\bigl(1-F(z)\bigr)^{1+\Delta} \bigr]\\[-1pt]
&=&1-2^{-\Delta}\vspace*{-1pt}
\end{eqnarray*}
and the maximum is attained at the point $z=F^{-1}(\frac12).$
The score generating function of the Psi test is $\varphi(u)=\ln u-\ln
(1-u),  0<u<1.$

Similarly, Savage \cite{Savage56}
proved that the \textit{Savage test} with the critical region
\[
\sum_{i=1}^n\sum_{j=R_{m+i}}^N\frac{1}{j}\leq C_{\alpha}\vspace*{-1pt}
\]
is the locally most powerful rank test of ${\mathbf H}_0$ against the
class of alternatives (\ref{alternatives}) with
%
\begin{equation}\label{Savage_alternative1}
G_{\Delta}^{(1)}(z)=F(z),  \qquad   G_{\Delta}^{(2)}(z)=F^{1+\Delta
}(z),  \qquad  z\geq0,\  \Delta>0.\vspace*{-1pt}
\end{equation}
Again, the $Y$'s are stochastically larger than the $X$'s, and the
Kolmogorov distance of~$G_{\Delta}^{(1)}$ and $G_{\Delta}^{(2)}$ is
equal to
\begin{eqnarray*}
d_K \bigl(G_{\Delta}^{(1)},G_{\Delta}^{(2)} \bigr)&=& \sup_{z\geq
0} \bigl[F(z)\bigl(1-F^{\Delta}(z)\bigr) \bigr]\\[-1pt]
&=& \Delta(1+\Delta)^{-1-{1}/{\Delta}},\vspace*{-1pt}
\end{eqnarray*}
whose maximum is attained at $z=F^{-1} ((1+\Delta)^{-1/\Delta
} ).$ The score generating function of this test is $\varphi
(u)=1+\ln u,  0<u<1.$

Assume that $F$ is increasing and let $U_{k}=F(z_{k}),  k=1,\ldots ,N.$
Under the alternati\-ve~(\ref{51}), the ranks $R_1,\ldots ,R_N$ are also
the ranks of the variables
$U_{1},\ldots ,U_{m},V_{m+1},\ldots ,V_{N},$
where $V_{k}=(1-\Delta) U_{k}+\Delta U^2_{k},  k=m+1,\ldots ,N.$ An
analogous consideration applies to the alternatives (\ref{G}) and
(\ref{Savage_alternative1}). Hence, the distribution of the ranks
$R_1,\ldots ,R_N$ is independent of $F$ (is distribution-free) under
the hypothesis as well as under the alternatives, and thus the power
functions of all rank tests against the alternatives (\ref{51}), (\ref
{G}) and~(\ref{Savage_alternative1}) are distribution-free.
The Lehmann alternatives can be well interpreted, are flexible and can
describe various experimental situations well. Besides the linear rank
tests, we can also consider the two-sample Kolmogorov--Smirnov test
based on the empirical distribution functions of the interpoint
distances, for the purposes of comparison. The randomized
Kolmogorov--Smirnov test, following a similar structure as the tests in
Section~\ref{sec2},
is also distribution-free.
Instead of interpoint distances, we can consider the rank tests based
on the depths using similar ideas.\vadjust{\goodbreak}

We shall concentrate on the two-sample Wilcoxon, Psi and Savage rank
tests because they are easy to perform, are locally most powerful and
are locally unbiased
against some alternatives of Lehmann type.
The ranks are distribution-free not only under the hypothesis, but also
under the Lehmann alternatives, hence the powers of the rank tests are
independent of the distribution of the data. This is an advantage 
because we do not need to calculate the distribution of the distances.
Several authors (e.g., \cite{FriedmanRafsky79,HenzePenrose99,Maa96,Rosenbaum05,Savage56})
considered various distances of two sets of multivariate observations
from some specified point, constructed the critical regions and verified
their consistencies against distant alternatives. However, the
questions of the finite-sample behavior of these tests, their
unbiasedness and against which alternatives, and their efficiency
against local alternatives, remains open.
If the test is not unbiased against some alternative of interest, then
its power can be less than the significance level, say less than
$\alpha=0.05,$ hence such a test is not suitable for verifying the
hypothesis against this specific alternative.

The sequences of alternatives (\ref{alternatives}) corresponding to
(\ref{Lehmann1}), (\ref{G}) and (\ref{Savage_alternative1})
are contiguous with respect to the sequence $\{\prod_{i=1}^N F(z_i)\}$
provided that $\Delta_N=N^{-1/2}\Delta_0$ with $\Delta_0$ fixed,
$0<\Delta_0<\infty,$ as shown in the \hyperref[appm]{Appendix}. Hence, we are able to
evaluate the local asymptotic powers of the tests; this is done in the
next section, along with the numerical illustration and comparison of
the tests.

\subsection{Local asymptotic powers of the tests} \label{section3}\label{sec3.1}

We shall assume throughout that
\[
\lim_{\Ny}\frac{m_N}{N}=\lambda\in(0,1).
\]
Let $U_i=F(Z_i),  i=1,\ldots ,N.$ The alternative (\ref
{alternatives}) in the special cases (\ref{51}), (\ref{G}) and~(\ref
{Savage_alternative1}) can then be rewritten as follows:
%
\begin{eqnarray}\label{55}
 \widetilde{G}_{\Delta}^{(1)}(u)&=&u, \qquad   \widetilde{G}_{\Delta
}^{(2)}(u)=(1-\Delta)u+\Delta u^2, \qquad   0\leq u\leq1,\nonumber\\
 \widetilde{G}_{\Delta}^{(1)}(u)&=&1 - (1 - u)^{1+\Delta}, \qquad
\widetilde{G}_{\Delta}^{(2)}(u)=\widetilde{G}(u,\Delta)=
u^{1+\Delta}, \nonumber
\\[-8pt]
\\[-8pt]
 \widetilde{G}_{\Delta}^{(1)}(u)&=&u, \qquad    \widetilde{G}_{\Delta
}^{(2)}(u)=\widetilde{G}(u,\Delta)=u^{1+\Delta},\nonumber\\
  \Delta&>&0, \qquad   0\leq u\leq1.
\nonumber
\end{eqnarray}
Because $\Delta$ is the parameter of interest and the alternatives
(\ref{55}) are contiguous with respect to the sequence of hypotheses
$ \{\prod_{i=1}^N u_iI[0\leq u_i\leq1] \}$ under $\Delta
_N=N^{-1/2}\Delta_0$ (see Appendix for the proof), we can study the
powers of the rank tests
under alternatives~(\ref{55}) without loss of generality.

Consider the centered test criterion
%
\begin{equation}\label{56}
S_N^*=N^{-1/2} \Biggl[-\frac{n}{N}\sum_{i=1}^m a_N(R_{Ni})+\frac
{m}{N}\sum_{i=m+1}^N a_N(R_{Ni}) \Biggr].\vadjust{\goodbreak}
\end{equation}
If the scores are generated by a nondecreasing score function $\varphi
\dvtx (0,1)\mapsto\mathbb{R}$ which is square-integrable on $(0,1),$
then the asymptotic distribution of (\ref{56})
under contiguous alternatives follows from the LeCam theorems (see
\cite{HajekSidak67} or \cite{HajekSidakSen99}).
Namely, $S_N^*$
will be asymptotically normally distributed ${\cal N}(\mu, \sigma^2)$ with
\begin{eqnarray*}
 \mu&=&\lambda(1-\lambda)\int_0^1\varphi(u)\varphi^*(u)\,\mathrm{d}u, \qquad
\varphi^*(u)=\frac{\partial\ln\tilde{g}(u,\Delta)}{\partial
\Delta}\bigg |_{\Delta=0},\\[-2pt]
 \tilde{g}(u,\Delta)&=&
\frac{\partial\tilde{G}(u,\Delta)}{\partial u}   \qquad \mbox{being
the density of }   \widetilde{G}(u,\Delta),\\[-2pt]
\sigma^2&=&\lambda(1-\lambda)\int_0^1\varphi^2(u)\,\mathrm{d}u.
\end{eqnarray*}

The test rejects ${\mathbf H}$ on the significance level $\alpha$ provided
$S_N^*\geq\sigma\Phi^{-1}(1-\alpha),$
where $\Phi$ is the standard normal distribution function. Hence, the
asymptotic power of the test under the alternative
${\mathbf K}_N\dvtx  \Delta_N=N^{-1/2}\Delta_0$ equals
\begin{eqnarray*}
&&\lim_{m,n\rightarrow\infty}{\Q}_{{\mathbf K}_N} \biggl(\frac
{S_N^*-\mu}{\sigma}\geq\Phi^{-1}(1-\alpha)-\frac{\mu}{\sigma
} \biggr)\\[-2pt] 
&& \quad =1-\Phi \biggl(\Phi^{-1}(1-\alpha)-\frac{\mu}{\sigma} \biggr)\\[-2pt]
&& \quad =1-\Phi \biggl(\Phi^{-1}(1-\alpha)-\frac{\sqrt{\lambda(1-\lambda
)}}{A}\Delta_0 \int_0^1\varphi(u)\varphi^*(u)\,\mathrm{d}u \biggr),
\end{eqnarray*}
where $A^2=\int_0^1\varphi^2(u)\,\mathrm{d}u.$
The relative asymptotic efficiency of a test $S_{N1}^*$ with respect to
a different test $S_{N2}^*$
is given as the ratio
%
\begin{equation}
\biggl (\frac{\mu^{(1)}}{\sigma_1}\bigg /\frac{\mu^{(2)}}{\sigma
_2} \biggr)^2,
\label{limitratio}
\end{equation}
where $\mu^{(1)}$ and $\sigma_1^2$ are, respectively, the asymptotic
mean and variance of the statistic~$S_{N1}^*$
and $\mu^{(2)}$ and $\sigma_2^2$ are those of $S_{N2}^*$.

\begin{table*}[b]
\vspace*{-3pt}
\tabcolsep=0pt
\tablewidth=285pt
\caption{Relative asymptotic efficiencies under various alternatives}
\label{tab1}\begin{tabular*}{285pt}{@{\extracolsep{\fill}}@{}llllll@{}}
\hline
Alternative &\multicolumn{5}{c@{}}{Test}\\[-5pt]
&\multicolumn{5}{c@{}}{\hrulefill}\\
& Wilcoxon & Psi & Savage & van der Waerden & Median\\
\hline
(\ref{51}) & 1.000 & 0.912 & 0.750 & 0.955 & 0.750 \\
(\ref{G}) & 0.912 & 1.000 & 0.882 & 0.992 & 0.584 \\
(\ref{Savage_alternative1}) & 0.750 & 0.822 & 1.000 & 0.816 & 0.480 \\
\hline
\end{tabular*}
\end{table*}

  Table~\ref{tab1}
summarizes the relative asymptotic efficiencies
of the Wilcoxon, Psi and Savage tests with respect to the locally
most
powerful rank test
for specified Lehmann alternatives. These values are computed with the\vadjust{\goodbreak}
aid of (\ref{limitratio}).
For the purposes of illustration, we also add the van der Waerden and
median tests, and their relative
asymptotic efficiencies with respect to the locally most powerful rank tests.

For the next illustration, consider the Lehmann alternative (\ref{51})
and compare the locally most powerful Wilcoxon test (the score function
$\varphi(u)=2u-1,  0\leq u\leq1$) with the Kolmogorov--Smirnov test.
The asymptotic power of the Wilcoxon test against ${\mathbf K}_N$ equals
%
\begin{eqnarray}\label{Wilcox_as_power1}
&&
\lim_{\min(m,n)\rightarrow\infty}
{\Q} \biggl(\biggl[\frac13\lambda(1-\lambda)\biggr]^{-1/2}
S_{N}\geq\Phi^{-1}(1-\alpha)\big |{\mathbf K}_N  \biggr)\nonumber
\\[-8pt]
\\[-8pt]
&& \quad
= 1-\Phi \Biggl(\Phi^{-1}(1-\alpha)-\Delta_0 \sqrt{\frac{\lambda
(1-\lambda)}{3}} \Biggr).
\nonumber
\end{eqnarray}

For small values $\Delta_0,$ it can be further approximated in the
following way:
%
\begin{eqnarray}\label{power_Wilcox_taylor1}
&&{\Q} \biggl(\biggl[\frac13\lambda(1-\lambda)\biggr]^{-1/2}S_{N}\geq
\Phi^{-1}(1-\alpha)\big |{\mathbf K}_N \biggr)\nonumber
\\[-8pt]
\\[-8pt]
&& \quad \approx \alpha+\Delta_0\cdot\Phi^{\top}\bigl(\Phi^{-1}(1-\alpha)\bigr)
\sqrt{\frac{\lambda(1-\lambda)}{3}}.
\nonumber
\end{eqnarray}
Let us now consider the Kolmogorov--Smirnov test against alternative
(\ref{51}).
Let $\hat{F}_{m}$ and $\hat{G}_{n}$ be the respective empirical
distribution functions of samples $X_1,\ldots ,X_m$ and
$Y_1,\ldots ,Y_n.$ Then, by H\'{a}jek \textit{et al.} \cite{HajekSidakSen99}, Theorem VI.3.2,
we have%
\begin{eqnarray*}
&&\lim_{m,n\rightarrow\infty}{\Q} \Biggl(
\sqrt{\frac{n   m}{n + m}} \sup_{x \in R}
 \bigl( \hat{G}_{n}(x) - \hat{F}_{m}(x)  \bigr) \geq
\sqrt{-\frac{1}{2} \log\alpha} |{\mathbf K_N} \Biggr) \nonumber
\\[-6.5pt]
\\[-6.5pt]
&& \quad ={\Q} \biggl(
\sup_{0 \leq u \leq1}
 \bigl(\mathcal{B}(u) + \Delta_0  \sqrt{\lambda( 1 - \lambda)}
u (1 - u)
 \bigr) \geq
\sqrt{-\frac{1}{2} \log\alpha} \biggr),
\nonumber
\end{eqnarray*}
where $\mathcal{B}(u)$ is a Brownian bridge. The last probability cannot easily be calculated analytically. Hence, we resort to a linear
approximation around the
point $\Delta_0 = 0$ and get
%
\begin{eqnarray} \label{power_KS_taylor1}
&&
{\Q} \biggl(
\sup_{0 \leq u \leq1}
 \bigl( \mathcal{B}(u) + \Delta_0  \sqrt{\lambda( 1 - \lambda)}
u (1 - u)
 \bigr) \geq
\sqrt{-\frac{1}{2} \log\alpha}
 \biggr) \nonumber
 \\[-8pt]
 \\[-8pt]
&& \quad
\approx
\alpha+2 \Delta_0  \sqrt{\lambda( 1 - \lambda)} \alpha
\sqrt{-\frac{1}{2} \log\alpha}
\int_{0}^{1}(2u - 1) \psi(\alpha, u)\,\mathrm{d}u,
\nonumber
\end{eqnarray}
where
\[
\psi(\alpha, u) =
2  \Phi \biggl(
\frac{(2u-1) \sqrt{- ({1}/{2}) \log\alpha}}{\sqrt{u (1-u)}}
 \biggr) - 1.
 \]
Table~\ref{tab2} gives the asymptotic powers (for $\alpha=0.05$) of the Wilcoxon
test (As.W) and the Kolmogorov--Smirnov test (As.KS) computed from
(\ref{Wilcox_as_power1}) and (\ref{power_KS_taylor1}); these powers
are compared with empirical powers (Obs.W, Obs.KS) obtained by simulations
of 30, 100, 500 and 1000 observations in both samples.
The simulations were carried out in the~{\textsf R} programming language
using 500\,000 replications
under the alternative (\ref{Lehmann1}), where~$F$ denotes the distribution
function of the uniform $\mathbf{R}(0,1)$ distribution. We recall that
the powers of
rank tests under Lehmann alternatives are also distribution-free for finite
samples.

\begin{table*}
\tabcolsep=0pt
\caption{Comparison of the empirical powers for various sample
sizes and of the local asymptotic
powers of the Wilcoxon and Kolmogorov--Smirnov tests against
the alternative (\protect\ref{Lehmann1})}
\label{tab2}\begin{tabular*}{\textwidth}{@{\extracolsep{\fill}}llllll@{\hspace{3pt}}lllll@{}}
\hline
$\Delta_{0}$ & \multicolumn{4}{c}{Obs.W, $m=n=$} & As.W &\multicolumn{4}{@{\hspace{3pt}}c}{Obs.KS, $m=n=$} & As.KS \\[-5pt]
& \multicolumn{4}{c}{\hrulefill} &  &\multicolumn{4}{c}{\hrulefill} &
 \\
& 30 & 100 & 500 & 1000 & & 30 & 100 & 500 & 1000\\
\hline
0.0 & 0.050 & 0.050 & 0.050 & 0.050 & 0.050 & 0.036 & 0.039 & 0.048 &
0.048 & 0.050 \\
0.1 & 0.052 & 0.053 & 0.053 & 0.053 & 0.053 & 0.038 & 0.040 & 0.050 &
0.052 & 0.053 \\
0.2 & 0.056 & 0.056 & 0.056 & 0.056 & 0.056 & 0.040 & 0.043 & 0.053 &
0.054 & 0.055 \\
0.3 & 0.059 & 0.059 & 0.059 & 0.059 & 0.060 & 0.042 & 0.046 & 0.056 &
0.057 & 0.058 \\
0.4 & 0.063 & 0.064 & 0.063 & 0.063 & 0.063 & 0.044 & 0.050 & 0.059 &
0.060 & 0.061 \\
0.5 & 0.067 & 0.066 & 0.067 & 0.066 & 0.067 & 0.047 & 0.052 & 0.063 &
0.063 & 0.063 \\
0.6 & 0.069 & 0.069 & 0.071 & 0.071 & 0.071 & 0.049 & 0.055 & 0.066 &
0.067 & 0.066 \\
0.7 & 0.076 & 0.074 & 0.074 & 0.075 & 0.075 & 0.053 & 0.057 & 0.069 &
0.070 & 0.069 \\
0.8 & 0.077 & 0.079 & 0.079 & 0.080 & 0.079 & 0.054 & 0.061 & 0.073 &
0.074 & 0.072 \\
0.9 & 0.081 & 0.083 & 0.083 & 0.082 & 0.083 & 0.057 & 0.063 & 0.076 &
0.077 & 0.074 \\
1.0 & 0.085 & 0.088 & 0.087 & 0.088 & 0.088 & 0.060 & 0.067 & 0.080 &
0.081 & 0.077 \\
2.0 & 0.141 & 0.142 & 0.141 & 0.143 & 0.143 & 0.100 & 0.107 & 0.126 &
0.131 & 0.104 \\
3.0 & 0.214 & 0.217 & 0.215 & 0.218 & 0.218 & 0.155 & 0.165 & 0.185 &
0.193 & 0.131 \\
\hline
\end{tabular*}
\vspace*{-3pt}
\end{table*}

The asymptotic approximation (\ref{Wilcox_as_power1}) of the power of
the Wilcoxon
test is already very good for $m=n=100$. Unfortunately, the linear
approximation (\ref{power_KS_taylor1})
of the power of the Kolmogorov--Smirnov test only works in a local
neighborhood of the null hypothesis as the power function increases
exponentially.
Even for small values of $\Delta_0,$ the approximation~(\ref{power_KS_taylor1}) of the power of the Kolmogorov--Smirnov test
is very good only
for large sample sizes.

Table~\ref{tab3} compares the slopes in
linear approximations of asymptotic powers of the Kolmogo\-rov--Smirnov
and Wilcoxon tests, given
in (\ref{power_Wilcox_taylor1}) and (\ref{power_KS_taylor1}), under
various sizes of the tests. The first column
gives the size of the test, the second column
the slope for the Kolmogorov--Smirnov test (K--S), the third column
gives the slope for
the Wilcoxon test and the last column gives the
ratio of the two slopes.

\begin{table*}
\tabcolsep=0pt
\tablewidth=7.6cm
\caption{Slopes of the Kolmogorov--Smirnov and Wilcoxon tests
at various levels of significance}
\label{tab3}
\begin{tabular*}{7.6cm}{@{\extracolsep{\fill}}llll@{}}
\hline
$\alpha$ & K--S & Wilcoxon & Wilcoxon/K--S \\
\hline
0.001 & 0.001 & 0.002 & 2.070 \\
0.010 & 0.009 & 0.015 & 1.680 \\
0.025 & 0.022 & 0.034 & 1.500 \\
0.050 & 0.044 & 0.059 & 1.350 \\
0.100 & 0.086 & 0.101 & 1.180 \\
\hline
\end{tabular*}
\end{table*}

%
\section{Numerical comparison of Hotelling- and Wilcoxon-type tests}\label{sec4}

The empirical powers of the Hotelling $T^2$ and Wilcoxon two-sample
tests are compared under bivariate\vadjust{\goodbreak} normal and Cauchy distributions with
various parameters; the Wilcoxon test of type (\ref{10}), (\ref{3})
is based on the ranks of the Euclidean interpoint distances. The
Hotelling test distinguishes well two normal samples contrasting
in locations, even if they also differ in scales. However, in some
situations, the Wilcoxon test even competes well with the Hotelling
test, namely, when either the samples differ only moderately in
locations or when they differ considerably in scales. This is
illustrated by Table~\ref{tab4}, which provides empirical powers of Hotelling
and Wilcoxon tests for a comparison of two bivariate normal samples.
The sample sizes are
$m=n=10,  100,  1000$ and the simulations are based on 10\,000 replications.
The first sample always has distribution ${\mathcal
N}_2({\boldgreek\mu}_1, {\boldgreek\Sigma}_1)$ with
${\boldgreek\mu}_1 = (0,0)^{\top}$ and ${\boldgreek\Sigma}_1 =
\diag\{1,1 \},$ while
the second sample has $\mathcal{N}_2({\boldgreek\mu}_2,\matice)$
with various parameters.

\begin{table*}
\tabcolsep=0pt
\caption{Powers of two-sample Hotelling $T^2$ test (H) and of
two-sample Wilcoxon test (W) based on distances for various $m=n$ and
$\alpha=0.05.$ The first sample always has ${\mathcal N}_2({\boldgreek
\mu}_1, {\boldgreek\Sigma}_1)$ distribution with
${\boldgreek\mu}_1 = (0,0)^{\top}$ and ${\boldgreek\Sigma}_1 =
\diag\{1,1 \}$.
The second sample has $\mathcal{N}_2({\boldgreek\mu}_2,\matice)$
with various ${\boldgreek\mu}_2,\matice$ specified in the first
column}
\label{tab4}
\vspace*{2pt}
\begin{tabular*}{\textwidth}{@{\extracolsep{\fill}}lllll@{}}
\hline
Second sample & Test & $m=n=10$ & $m=n=100$ & $m=n=1000$\\
\hline
${\boldgreek\mu}_2= (0,0)^T$ & H & 0.0471 & 0.0481 & 0.0493 \\
$\matice=\diag\{1, 1\}$ & W & 0.0457 & 0.0505 & 0.0487 \\
[6pt]
${\boldgreek\mu}_2 =(0.2, 0.2)^T$ & H & 0.0771 & 0.4115 & 1.0000 \\
$\matice=\diag\{1, 1\}$& W & 0.0520 & 0.1715 & 0.6458 \\
${\boldgreek\mu}_2 =(0.5, 0.5)^T$ & H & 0.2318 & 0.9962 & 1.0000 \\
$\matice=\diag\{1, 1\}$ & W & 0.1085 & 0.5701 & 0.8617 \\
[6pt]
${\boldgreek\mu}_2 =(0, 0)^T$ & H & 0.0659 & 0.0561 & 0.0452 \\
$\matice=\diag\{0.1, 0.1\}$ & W & 0.7994 & 0.9998 & 1.0000 \\
${\boldgreek\mu}_2 =(0, 0)^T$ & H & 0.0653 & 0.0456 & 0.0530 \\
$\matice=\diag\{0.2, 0.2\}$ & W & 0.4851 & 0.9932 & 1.0000 \\
${\boldgreek\mu}_2=(0, 0)^T$ & H & 0.0521 & 0.0521 & 0.0463 \\
$\matice=\diag\{0.5, 0.5\}$ & W & 0.1182 & 0.7034 & 0.9968 \\
${\boldgreek\mu}_2 =(0, 0)^T$ & H & 0.0531 & 0.0530 & 0.0514 \\
$\matice=\diag\{1.5, 1.5\}$ & W & 0.0656 & 0.2881 & 0.8525 \\
${\boldgreek\mu}_2=(0, 0)^T$ & H & 0.0552 & 0.0518 & 0.0508 \\
$\matice=\diag\{2, 2\}$& W & 0.0999 & 0.5395 & 0.9670 \\
${\boldgreek\mu}_2=(0,0)^T$ & H & 0.0572 & 0.0546 & 0.0521\\
$\matice=\diag\{1.0, 0.2\}$ & W & 0.1029 & 0.6568 & 0.9936\\
[6pt]
${\boldgreek\mu}_2=(0.1, 0.1)^T$ & H & 0.0553 & 0.1266 & 0.7897 \\
$\matice=\diag\{1.1, 1.1\}$ & W & 0.0491 & 0.0932 & 0.4232 \\
${\boldgreek\mu}_2 =(0.1, 0.1)^T$ & H & 0.0601 & 0.1167 & 0.7183 \\
$\matice=\diag\{1.5, 1.5\}$ & W & 0.0667 & 0.3182 & 0.7690 \\
${\boldgreek\mu}_2 =(0.2,0.2)^T$ & H & 0.0742 & 0.3656 & 1.0000\\
$\sigma_1^2=1, \sigma_2^2=1.5$ & W & 0.0548 & 0.2246 & 0.6907\\
$ {\boldgreek\mu}_2=(0.2, 0.2)^T$ & H & 0.0710 & 0.3402 & 0.9994 \\
$\matice=\diag\{1.5, 1.5\}$ & W & 0.0668 & 0.3551 & 0.7597 \\
\hline
\end{tabular*}
\vspace*{12pt}
\end{table*}

\begin{table*}
\tabcolsep=0pt
\caption{Powers of two-sample Hotelling $T^2$ test (H) and
two-sample Wilcoxon test (W) based on distances for various $m=n$ and
$\alpha=0.05.$
The first sample $\X$ always has the two-dimensional Cauchy
distribution. The second sample $\mathbf{Y}$ is obtained as $\mathbf
{Y} = \mm+ \sigma\Y^*$, where $\Y^*$ is generated as a
two-dimensional Cauchy distribution independent of $\X.$ Values of
$m,n,\mm$ and $\sigma$ are specified in the first column}
\label{tab5}\begin{tabular*}{\textwidth}{@{\extracolsep{\fill}}llllll@{}}
\hline
Second sample & Test & $m=n=10$ & $m=n=25$ & $m=n=100$ & $m=n=1000$\\
\hline
$\mm= (0,0)^T$ & H & 0.0191 & 0.0156 & 0.0171 & 0.0217 \\
$\sigma=1$ & W & 0.0450 & 0.0478 & 0.0510 & 0.0442\\
[3pt]
$\mm=(0.2, 0.2)^T$ & H & 0.0227 & 0.0232 & 0.0227 & 0.0174 \\
$\sigma=1$ & W & 0.0468 & 0.0536 & 0.0874 & 0.3925\\
$\mm=(0.5, 0.5)^T$ & H & 0.0408 & 0.0404 & 0.0414 & 0.0361\\
$\sigma=1$ & W & 0.0664 & 0.1115 & 0.2937 & 0.7470\\
$\mm=(1, 1)^T$ & H & 0.1038 & 0.1193 & 0.1260 & 0.1226 \\
$\sigma=1$ & W & 0.1219 & 0.2710 & 0.6235 & 0.8893 \\
$\mm=(5, 5)^T$ & H & 0.7387 & 0.7535 & 0.7683 & 0.7772 \\
$\sigma=1$ & W & 0.7574 & 0.9441 & 0.9782 & 0.9944 \\
[3pt]
$\mm= (0,0)^T$ & H & 0.0200 & 0.0171 & 0.0193 & 0.0103\\
$\sigma=1.5$ & W & 0.0664 & 0.1207 & 0.3419 & 0.8428 \\
$\mm= (0,0)^T$ & H & 0.0207 & 0.0168 & 0.0182 & 0.0172 \\
$\sigma=2$ & W & 0.1082 & 0.2439 & 0.6123 & 0.9135 \\
[3pt]
$\mm= (0.2, 0.2)^T$ & H & 0.0189 & 0.0201 & 0.0196 & 0.0240\\
$\sigma=1.5$ & W & 0.0710 & 0.1297 & 0.3495 & 0.8249 \\
$\mm= (1, 1)^T$ & H & 0.0741 & 0.0814 & 0.0865 & 0.0943\\
$\sigma=1.5$ & W & 0.1088 & 0.2188 & 0.4925 & 0.8462 \\
$\mm= (2, 2)^T$ & H & 0.2356 & 0.2546& 0.2690 & 0.2716\\
$\sigma=1.5$ & W & 0.2092 & 0.4139 & 0.7395 & 0.9259 \\
[3pt]
$\mm= (0.2, 0.2)^T$ & H & 0.0248 & 0.0186 & 0.0217 & 0.0158\\
$\sigma=2$ & W & 0.1134 & 0.2401 & 0.5990 & 0.9151 \\
$\mm= (1, 1)^T$ & H & 0.0575 & 0.0616 & 0.0623 & 0.0676\\
$\sigma=2$ & W & 0.1330 & 0.2797 & 0.5619 & 0.8272 \\
$\mm= (2, 2)^T$ & H & 0.1796 & 0.1936& 0.2045 & 0.2164\\
$\sigma=2$ & W & 0.1771 & 0.3513 & 0.6531 & 0.8981 \\
\hline
\end{tabular*}\vspace*{2pt}
\end{table*}

We also refer to the simulation study of   \cite{ZuoHe06}
which compared the empirical powers of the Liu--Singh rank-sum test ($Q$)
based on the depths, the Hotelling and the Hetmansperger \textit{et al.} \cite
{HetOja98} tests for two bivariate normal samples. Under normality, the
$Q$-test mostly dominates the other two tests, as well as the Wilcoxon
test based on interpoint distances.
However, the (local) unbiasedness of the $Q$-test against two-sample
alternatives is doubtful under asymmetric distributions of the depths,
while a one-sided alternative in depths has a difficult interpretation
in the original data.

Table~\ref{tab5} presents the empirical powers of the tests comparing two
samples from the bivariate Cauchy distributions. The first sample
$\mathbf{X}$ has a two-dimensional Cauchy distribution with
independent components. The second sample $\mathbf{Y}$ is obtained as
a random sample ${\mathbf Y^*}$ from the two-dimensional Cauchy
distribution with independent components, independent of $\mathbf{X}$,
transformed to $\mathbf{Y} = \mm+ \sigma\Y^*$ for certain shifts
$\mm$ and scales $\sigma.$ The results are based on 10\,000
replications. The Wilcoxon test is far more powerful than the Hotelling
test, already under a small shift. The Hotelling test fails completely
if $\mm=0$ but $\sigma\neq1,$ while the Wilcoxon test still
distinguishes well the samples. The Wilcoxon test also dominates the
Hotelling test in other situations.

The rank tests based on interpoint distances are distribution-free,
both under the hypothesis and under the Lehmann alternatives, while the
exact distribution of the distances can remain unknown when performing
the tests. The tests are locally unbiased against one-sample
alternatives. If the interpoint distances are replaced with other
scalar characteristics which are symmetrically distributed under the
hypothesis, then the tests are also locally unbiased against the
two-sample alternatives. The Lehmann alternatives reflect the practical
situations well.

\begin{appendix}
\section*{Appendix: Contiguity of Lehmann's alternatives}\vspace*{5pt}\label{appm}

Let $\{P_{N1},\ldots ,P_{NN}\}$ and $\{Q_{N1},\ldots ,Q_{NN}\}$ be two
triangular arrays of probability measures
defined on the measurable space $({\mathcal X},{\mathcal A}),$ and let
$P_N^{(N)}=\prod_{i=1}^N P_{Ni}$ and
$Q_N^{(N)}=\prod_{i=1}^N Q_{Ni}$ denote the respective product
measures, $N=1,2,\ldots .$ Further, denote by $p_{Ni}$ and $q_{Ni}$
the respective densities of $P_{Ni}$ and $Q_{Ni}$ with respect to a
$\sigma$-finite measure $\mu_i$, which can also be
$\mu_i=P_{Ni}+Q_{Ni},  i=1,\ldots ,N.$

Oosterhoff and van Zwet \cite{Oosterhoff79}
proved that $\{Q_N^{(N)}\}$ is contiguous with respect to $\{P_N^{(N)}\}$
if and only if
%
\begin{equation}\label{Zwet1}
\limsup_{N\rightarrow\infty}\sum_{k=1}^NH^2(P_{Nk},Q_{Nk})<\infty
\end{equation}
and
%
\begin{equation}\label{Zwet2}
\lim_{N\rightarrow\infty}\sum_{k=1}^N Q_{Nk} \biggl\{\frac
{q_{Nk}(X_{Nk})}{p_{Nk}(X_{Nk})}\geq c_N \biggr\}=0 \qquad
\forall c_N\rightarrow\infty,
\end{equation}
where
\[
H(P,Q)= \biggl[\int \bigl(\sqrt{p}-\sqrt{q} \bigr)^2\,\mathrm{d}\mu
\biggr]^{1/2}= \biggl[2\int \bigl(1-\sqrt{pq} \bigr)\,\mathrm{d}\mu \biggr]^{1/2}
\]
is the Hellinger distance of $P,Q.$

Put $\Delta_N=N^{-1/2}\Delta_0$ with $\Delta_0>0$ fixed.
Applying (\ref{Zwet1}) and (\ref{Zwet2}), we can verify the
contiguity of the sequence
$ \{\prod_{k=1}^m G_{\Delta_N}^{(1)}(z_k)\prod_{k=m+1}^N
G_{\Delta_N}^{(2)}(z_k) \}$ with respect
to $ \{\prod_{k=1}^N F(z_k) \}$ for the alternatives (\ref
{Lehmann1}), (\ref{G}) and (\ref{Savage_alternative1}).
\begin{LEM}\label{Lemma1}
\textup{(i)} Let
%
\begin{equation}\label{a}
 \Biggl\{\prod_{k=1}^N F(z_k) \Biggr\}_{N=1}^{\infty}, \qquad   z_1,\ldots
,z_N\geq0,\   N=1,2,\ldots ,
\end{equation}
and
%
\begin{equation} \label{b}
\Biggl \{\prod_{k=1}^m G_{\Delta_N}^{(1)}(z_k)\prod_{k=m+1}^N
G_{\Delta_N}^{(2)}(z_k) \Biggr\}_{N=1}^{\infty},  \qquad  z_1,\ldots
,z_N\geq0,\  N=1,2,\ldots,
\end{equation}
be two sequences of probability distributions satisfying
\begin{eqnarray*}
 \Delta_N&=&N^{-1/2}\Delta_0>0, \qquad   \lim_{\Ny}\min\{m,n\}=\infty
,\\
 \lim_{\Ny}\frac mN&=&\lim_{\Ny}\frac{m_N}{N}=\lambda\in(0,1),
\end{eqnarray*}
where $G_{\Delta}^{(1)},  G_{\Delta}^{(2)}$ are given by either
(\ref{51}), (\ref{G}) or (\ref{Savage_alternative1}).
The sequence (\textup{\ref{b}}) is then contiguous with respect to the sequence
(\textup{\ref{a}}).
\end{LEM}

\begin{pf}
 (i) Let us first consider the Lehmann alternatives
(\ref{51}). Then,
\begin{eqnarray*}
&&\sum_{k=m+1}^N H^2(F(z_k),G_{\Delta_N}(z_k))\\
&& \quad =n\cdot\int_0^{\infty} f(z) \bigl[\sqrt{1+\Delta
_N\bigl(2F(z)-1\bigr)}-1 \bigr]^2\,\mathrm{d}z\\
&& \quad =n\int_0^{\infty} f(z)\frac{[1+\Delta_N(2F(z)-1)-1]^2}{
[\sqrt{1+\Delta_N(2F(z)-1)}+1 ]^2}\,\mathrm{d}z\\
&& \quad \leq n\Delta_N^2\int_0^{\infty}f(z)\bigl(2F(z)-1\bigr)^2\,\mathrm{d}z\\
&& \quad =4n\Delta_N^2\int_0^1\biggl(u-\frac12\biggr)^2=\frac13 n\Delta_N^2=
\frac13 \lambda_N\Delta_0<\infty
\end{eqnarray*}
and
\begin{eqnarray*}
&&\lim_{N\rightarrow\infty}\sum_{k=1}^N Q_{Nk}\biggl \{\frac
{q_{Nk}(X_{Nk})}{p_{Nk}(X_{Nk})}\geq c_N \biggr\}\\
&& \quad =0+\lim_{N\rightarrow\infty}\sum_{k=m+1}^N\int_0^{\infty}
I\bigl [1+\Delta_N\bigl(2F(z_k)-1\bigr)\geq c_N \bigr]\\
&&  \quad \hphantom{=0+\lim_{N\rightarrow\infty}\sum_{k=m+1}^N\int_0^{\infty}}    {}\times\bigl[1+\Delta
_N\bigl(2F(z_k)-1\bigr)\bigr]f(z_k)\,\mathrm{d}z_k=0
\end{eqnarray*}
because $c_N>1+N^{-1/2}\Delta_0$ for $n>N_0$ whenever $c_N\rightarrow
\infty.$ The contiguity is
thus verified.

(ii) Let
\[
G_{\Delta_N}(z_i)=
\cases{\displaystyle
1 - \bigl(1 - F(z_i)\bigr)^{1+\Delta_{N}}, & \quad $i\leq m$,\cr\displaystyle
(F(z_i))^{1+\Delta_{N}}, & \quad  $i\geq m+1$.
}
\]
Then,
\begin{eqnarray*}
&&\sum_{i=1}^N H^2(F(z_i),G_{\Delta_N}(z_i))
    \\[-2pt]
&& \quad =m\cdot\int_0^{\infty} f(z)\bigl [\sqrt{(1+\Delta
_N)\bigl(1-F(z)\bigr)^{\Delta_N}}-1 \bigr]^2\,\mathrm{d}z\\[-2pt]
&& \qquad{} +n\cdot\int_0^{\infty} f(z)\bigl [\sqrt{(1+\Delta
_N)(F(z))^{\Delta_N}}-1 \bigr]^2\,\mathrm{d}z\\[-2pt]
&& \quad \leq m\cdot\int_0^{\infty} f(z) \bigl[(1+\Delta_N)\bigl(1-F(z)\bigr)^{\Delta
_N}-1 \bigr]^2\,\mathrm{d}z\\[-2pt]
&& \qquad{} +n\cdot\int_0^{\infty} f(z) [(1+\Delta_N)(F(z))^{\Delta
_N}-1 ]^2\,\mathrm{d}z\\[-2pt]
&& \quad =m\cdot\int_0^1 [(1+\Delta_N)(1-u)^{\Delta_N}-1 ]^2\,\mathrm{d}z+
n\cdot\int_0^1 [(1+\Delta_N)u^{\Delta_N}-1 ]^2\,\mathrm{d}z\\[-2pt]
&& \quad =N\cdot\int_0^1 [(1+\Delta_N)u^{\Delta_N}-1 ]^2\,\mathrm{d}z\leq
\Delta_N<\infty
\end{eqnarray*}
and hence (\ref{Zwet1}) is proved for the alternative (\ref{G}).
Concerning (\ref{Zwet2}), we have
\begin{eqnarray*}
&&\lim_{N\rightarrow\infty}\sum_{i=1}^N Q_{Ni}\biggl \{\frac
{q_{Ni}(X_{Ni})}{p_{Ni}(X_{Ni})}\geq c_N \biggr\}\\
&& \quad =\lim_{\Ny} \biggl\{m\cdot\int_0^{\infty}I \bigl[(1+\Delta
_N)\bigl(1-F(z)\bigr)^{\Delta_N}\geq c_N \bigr]
(1+\Delta_N)\bigl(1-F(z)\bigr)^{\Delta_N}f(z)\,\mathrm{d}z\\
&& \quad \hphantom{=\lim_{\Ny} \biggl\{}{}+n\cdot\int_0^{\infty}I [(1+\Delta_N)(F(z))^{\Delta_N}\geq
c_N ](1+\Delta_N)(F(z))^{\Delta_N}f(z)\,\mathrm{d}z \biggr\}\\
&& \quad =\lim_{\Ny}N\cdot\int_0^1I [(1+\Delta_N)u^{\Delta_N}\geq
c_N ](1+\Delta_N)u^{\Delta_N}\,\mathrm{d}u\\
&& \quad =\lim_{\Ny}N\cdot\int_0^{1+\Delta_N}I [v\geq c_N
]v^{1/\Delta_N}\Delta_N^{-1}(1+\Delta_N)^{-1/\Delta_N}\,\mathrm{d}v=0
\end{eqnarray*}
because the set $\{v\dvtx   c_N\leq v\leq1+\Delta_N\}$ is empty for
$N>N_0.$

(iii) Similarly, for the alternative (\ref{Savage_alternative1}),
we have
\begin{eqnarray*}
&&
\sum_{k=m+1}^N H^2(F(z_k),G_{\Delta_N}(z_k))\\
&& \quad = n \cdot\int_0^{\infty} f(z)  \bigl[\sqrt{(1+ \Delta_N) F^{\Delta
_N}(z)} -1  \bigr]^2\,\mathrm{d}z\\
&& \quad
= n \int_0^{\infty}
f(z)\frac{[(1- \Delta_N) F^{\Delta_N}(z) - 1]^2}{ [\sqrt{(1+
\Delta_N)
F^{\Delta_N}(z)}+1 ]^2}\,\mathrm{d}z\\
&& \quad
\leq n \int_0^{1} [(1+\Delta_N) u^{\Delta_N} - 1 ]^2\,\mathrm{d}u
\leq
n \int_0^{1}  [\Delta_N^{2}+ (u^{\Delta_N} - 1)^{2}  ]\,\mathrm{d}u
\\
&& \quad
\leq
n  \biggl\{
\Delta_N^{2} +  \biggl[
\frac{u^{1+2 \Delta_N}}{1+2 \Delta_N} -
\frac{2 u^{1+\Delta_N}}{1+\Delta_N} + u
 \biggr]_{0}^{1}
 \biggr\}
=
n  \biggl\{
\Delta_N^{2} + \frac{2   \Delta_N^{2}}{(1 + 2 \Delta_N)(1 +
\Delta_{n})}
 \biggr\} \\
&& \quad
\leq
7  n   \Delta_N^{2} = 7  \lambda_{N}   \Delta_{0}^{2} < \infty.
\end{eqnarray*}
Condition~\eqref{Zwet2}
is verified analogously as for the alternative (\ref{51}).
\end{pf}
\end{appendix}

\section*{Acknowledgments}
The authors would like to thank the Editor and two Referees for their
valuable comments which helped to provide a better understanding of the
whole text. They also wish to thank Pranab K. Sen, Hannu Oja and Marek
Omelka for valuable discussions and M. Omelka also
for his help with calculating the powers of some tests.
This research was supported by the Project LC06024 of Ministry of Education, Youth and Sports of Czech Republic.
     J. Jure\v{{c}}kov\'{a} was also supported by the grant IAA101120801
of the Academy of Science of Czech Republic, by the Czech Republic Grant 201/09/0133 and by the research project
MSM 0021620839 of the Ministry of Education, Youth and Sports of Czech Republic.

\printhistory

\end{document}